\documentclass[11pt,a4paper]{article}

\usepackage{amsmath,amssymb,latexsym}
\usepackage{graphicx}
\usepackage{epsfig}
\usepackage{longtable}
\usepackage{amsthm}
\usepackage{xcolor}
\usepackage{textcomp}

\newtheorem{theorem}{Theorem}

\newtheorem{example}{Example}

\begin{document}

\begin{center}
{\Large\bf Extreme values of the mass distribution associated with a tetravariate quasi-copula}
\end{center}
\smallskip

\begin{center}
{\large\bf Manuel \'Ubeda-Flores
}
\end{center}

\begin{center}
{\it Department of Mathematics, University of Almer\'{\i}a, 04120 Almer\'{\i}a, Spain.\\
{\rm {\rm mubeda@ual.es}}}
\end{center}
\smallskip

\begin{center}{\bf Abstract}\end{center}

\noindent
In this note we study the extremes of the mass distribution associated with a tetravariate quasi-copula and compare our results with the
bi- and trivariate cases, showing the important differences between them.
\bigskip

\medskip


\section{Introduction}

$\phantom{999}$Alsina {\it et al.} \cite{Alsina1993} introduced the notion of a {\it quasi-copula} ---a more general concept than that of a copula \cite{Durante2016,Nelsen2006}--- in order to characterize operations on distribution functions that can, or cannot, be derived from operations on random variables defined on the same probability space. In the last few years these functions have attracted an increasing interest by researchers in some topics of fuzzy sets theory, such as preference modeling, similarities and fuzzy logics. For a complete overview of quasi-copulas, we refer to \cite{Arias2020,Sempi2017}.

Cuculescu and Theodorescu \cite{Cuculescu2001} characterize a multivariate quasi-copula, or $n$-quasi-copula, ---a 2-quasi-copula was previously characterized in \cite{Genest1999}--- as a function $Q\colon [0,1]^{n}\to [0,1]$ that satisfies:
\begin{enumerate}
\item  {\it boundary conditions}: $Q(u_{1},\ldots,u_{i-1},0,u_{i+1},\ldots,u_{n})=0$ and $Q(1,\ldots,1,u_i,1,\ldots,1)=u_{i}$, for any $(u_1,\ldots,u_n)\in [0,1]^n$;
\item  {\it monotonicity}: $Q$ is increasing in each variable; and
\item  {\it Lipschitz condition}: $\vert Q(u_1,\ldots,u_n)-Q(v_1,\ldots,v_n)\vert \le \sum_{i=1}^{n}\vert u_{i}-v_{i}\vert$, for any $(u_1,\ldots,u_n)$ and $(v_1,\ldots,v_n)$
in $[0,1]^{n}$.
\end{enumerate}

Every $n$-quasi-copula $Q$ satisfies the following inequalities: $W^{n}(u_1,\ldots,u_n)= \max(\sum_{i=1}^{n}u_i-n+1,0)\le Q(u_1,\ldots,u_n)\leq \min(u_1,\ldots,u_n)$. $W^{2}$ is a 2-copula and $W^n$, $n\ge 3$, is a {\it proper} $n$-quasi-copula, i.e., an $n$-quasi-copula but not an $n$-copula.

Consider an $n$-quasi-copula $Q$ and an $n$-{\it box} $B=[a_1,b_1]\times\cdots\times[a_n,b_n]$ in $[0,1]^{n}$. The $Q$-{\it volume} of $B$ is defined as $V_{Q}(B)=\sum
\mathrm{sgn}(c_1,\ldots,c_n)\,Q(c_1,\ldots,c_n),$ where the sum is taken over all the {\it vertices} $(c_{1},\ldots,c_{n})$ of $B$ ---i.e., each $c_k$ is equal to either $a_k$ or $b_k$--- and $\mathrm{sgn}(c_1,\ldots,c_n)$ is $1$ if $c_{k}=a_{k}$ for an even number of $k'$s, and $-1$ if $c_{k}=a_{k}$ for an odd number of $k'$s. We refer to $V_Q$ as the {\em mass distribution} associated with $Q$ (on $n$-boxes), and to $V_Q(B)$ as the mass accumulated by $Q$ on $B$.

One of the most important occurrences of quasi-copulas in statistics is due to the following observation \cite{Nelsen2004,RodUb2004}: Every set ${\cal{S}}$ of (quasi-)copulas has the smallest upper bound and the greatest lower bound in the set of quasi-copulas (in the sense of pointwisely ordered functions). These bounds do not necessarily belong to the set ${\cal{S}}$, nor they are necessarily copulas if the set consists of copulas only.

Differences between copulas and quasi-copulas are studied \cite{Fernandez2014,Fernandez2011,Nelsen2010}, and to understand the behavior of an $n$-quasi-copula, the step from studying a certain property of a quasi-copula from dimension 2 to 3 is crucial, since, in most of the cases known in the literature, although the behavior of the bivariate case to the trivariate case differs, this is not usually the case for the three-dimensional case and higher.

The mass distribution is a popular object of investigation \cite{DeBaets2007,Fernandez2023,Nelsen2005,Nelsen2002,Omladic2022,RodUb2009,Sempi2023}. In the case of investigating the extremes of the mass distribution of an $n$-quasi-copula, although the bivariate and trivariant cases differ, it might be expected that this would not be the case for the tri- and tetravariate dimensions; however, we will show that we obtain results that reveal this difference, which is, in principle, surprising, since the values of these extremes exceed a priori reasonable bounds. In the case $n=2$, there exists a unique $2$-box on which the minimal mass (which turns out to be $-1/3$) can be accumulated, as well as a unique $2$-box (the unit square itself) on which the maximal mass $1$ can be accumulated \cite{Nelsen2002}. In the case $n=3$ there still exists a unique $3$-box on which the minimal mass (which now turns out to be $-4/5$) can be accumulated, while there exist multiple $3$-boxes on which the maximal mass $1$ can be accumulated \cite{DeBaets2007}. The proof given in \cite{Nelsen2002} for the case $n=2$ is direct, while the one given in \cite{DeBaets2007} for the case $n=3$ is given by formulating this optimization problem as a linear programming problem.

In principle, the linear programming methodology can be applied to the case $n>3$ as well, with an obvious increase in complexity of the formulation, and, as we will show,
several surprises will appear (Section \ref{sec:main}). It is probable that other techniques can prove the result for any $n>3$ (for example, using characterizations of the positive volumes of an $n$-quasi-copula shown in \cite{Omladic2022,RodUb2009}), but no progress has been made so far in this regard. As a consequence of the results obtained, in Section \ref{sec:discus} we will make a conjecture about the extremes of the mass distribution of an $n$-quasi-copula for any $n>4$.

\section{The mass distribution associated with a 4-quasi-copula}\label{sec:main}

$\phantom{999}$In this section we study the extremes of the mass distribution associated with a 4-quasi-copula. We provide the main result of this note.

\begin{theorem}\label{th:main}Let $Q$ be a $4$-quasi-copula and let $B$ be a $4$-box in $[0,1]^4$. Then it holds $-9/7\le V_Q(B)\le 2$.
\end{theorem}

\begin{proof}Consider the 4-box $B=\left[a_1,b_1\right]\times \left[a_2,b_2\right]\times\left[a_3,b_3\right]\times\left[a_4,b_4\right]$, and let $Q$ be a 4-quasi-copula. In order to simplify the notation, we define by the following parameters the length of each of the edges of the 4-box: $a=b_1-a_1$, $b=b_2-a_2$, $c=b_3-a_3$ and $d=b_4-a_4$. The values of the 4-quasi-copula at the vertices of the 4-box $B$ are: $Q(a_1,a_2,a_3,a_4)=e$, $Q(a_1,a_2,a_3,b_4)=f+e$, $Q(a_1,a_2,b_3,a_4)=g+e$, $Q(a_1,a_2,b_3,b_4)=h+e$, $Q(a_1,b_2,a_3,a_4)=i+e$, $Q(a_1,b_2,a_3,b_4)=j+e$, $Q(a_1,b_2,b_3,a_4)=k+e$, $Q(a_1,b_2,b_3,b_4)=l+e$, $Q(b_1,a_2,a_3,a_4)=m+e$, $Q(b_1,a_2,a_3,b_4)=n+e$, $Q(b_1,a_2,b_3,a_4)=o+e$, $Q(b_1,a_2,b_3,b_4)=p+e$, $Q(b_1,b_2,a_3,a_4)=q+e$, $Q(b_1,b_2,a_3,b_4)=r+e$, $Q(b_1,b_2,b_3,a_4)=s+e$, $Q(b_1,b_2,b_3,b_4)=t+e$. We consider 24 fundamental parameters, namely: $a_1, a_2, a_3, a_4, a, b, c, d, e, f, g, h, i, j, k, l, m, n, o, p, q, r, s, t$. These parameters, along with the values $a+a_1$, $b+a_2$, $c+a_3$ and $d+a_4$, and applying the boundary conditions of a 4-quasi-copula, must lie in the interval $[0,1]$.

We apply the other conditions of a 4-quasi-copula. Since $Q$ is non-decreasing in each variable and satisfies the 1-Lipschitz condition, we have\smallskip

$m\le a$,\quad $i\le b$,\quad $g\le c$,\, and\, $f\le d,$\smallskip

\noindent
and\smallskip

$\max(g,f)\le h\le \min(d+g,c+f)$,

$\max(i,f)\le j\le \min(b+f,d+i)$,

$\max(i,g)\le k\le \min(b+g,c+i)$,

$\max(f,m)\le n\le \min(a+f,d+m)$,

$\max(g,m)\le o\le \min(a+g,c+m)$,

$\max(i,m)\le q\le \min(a+i,b+m)$,

$\max(h,j,k)\le l\le \min(b+h,c+j,d+k)$,

$\max(h,n,o)\le p\le \min(a+h,c+n,d+o)$,

$\max(j,n,q)\le r\le \min(a+j,b+n,d+q)$,

$\max(k,o,q)\le s\le \min(a+k,b+o,c+q)$,

$\max(l,p,r,s)\le t\le \min(a+l,b+p,c+r,d+s)$.\smallskip

Since $\max(u+v+w+z-3,0)\le Q(u,v,w,z)\le \min(u,v,w,z)$ for all $(u,v,w,z)\in[0,1]^4$, then, for each of the 16 vertices of the 4-box $B$ we have:\smallskip

$a_1+a_2+a_3+a_4-3\le e\le \min\left(a_1,a_2,a_3,a_4\right)$,

$a_1+a_2+a_3+a_4+d-3\le f+e\le \min\left(a_1,a_2,a_3,a_4+d\right)$,

$a_1+a_2+a_3+a_4+c-3\le g+e\le \min\left(a_1,a_2,a_3+c,a_4\right)$,

$a_1+a_2+a_3+a_4+c+d-3\le h+e\le \min\left(a_1,a_2,a_3+c,a_4+d\right)$,

$a_1+a_2+a_3+a_4+b-3\le i+e\le \min\left(a_1,a_2+b,a_3,a_4\right)$,

$a_1+a_2+a_3+a_4+b+d-3\le j+e\le \min\left(a_1,a_2+b,a_3,a_4+d\right)$,

$a_1+a_2+a_3+a_4+b+c-3\le k+e\le \min\left(a_1,a_2+b,a_3+c,a_4\right)$,

$a_1+a_2+a_3+a_4+b+c+d-3\le l+e\le \min\left(a_1,a_2+b,a_3+c,a_4+d\right)$,

$a_1+a_2+a_3+a_4+a-3\le m+e\le \min\left(a_1+a,a_2,a_3,a_4\right)$,

$a_1+a_2+a_3+a_4+a+d-3\le n+e\le \min\left(a_1+a,a_2,a_3,a_4+d\right)$,

$a_1+a_2+a_3+a_4+a+c-3\le o+e\le \min\left(a_1+a,a_2,a_3+c,a_4\right)$,

$a_1+a_2+a_3+a_4+a+c+d-3\le p+e\le \min\left(a_1+a,a_2,a_3+c,a_4+d\right)$,

$a_1+a_2+a_3+a_4+a+b-3\le q+e\le \min\left(a_1+a,a_2+b,a_3,a_4\right)$,

$a_1+a_2+a_3+a_4+a+b+d-3\le r+e\le \min\left(a_1+a,a_2+b,a_3,a_4+d\right)$,

$a_1+a_2+a_3+a_4+a+b+c-3\le s+e\le \min\left(a_1+a,a_2+b,a_3+c,a_4\right)$,

$a_1+a_2+a_3+a_4+a+b+c+d-3\le t+e\le \min\left(a_1+a,a_2+b,a_3+c,a_4+d\right)$.\smallskip

Now we apply linear programming considering to minimize the objective function
$$V_Q(B)=h+j+k+n+o+q+t-f-g-i-l-m-p-r-s$$
for a total of 192 linear inequality constraints, i.e., looking for the 4-boxes $B$ that minimize $V_Q(B)$. The simplex method ---see, e.g., \cite{Chong2001}--- shows the solution $V_Q(B)=-9/7$, and where we have used the Mathematica\textregistered\, package.

Analogously, by maximizing $V_Q(B)$, or equivalently, by minimizing $-V_Q(B)$, we obtain $V_Q(B)=2$, which completes the proof.
\end{proof}

The linear programming in the proof of Theorem \ref{th:main} shows the following solution for the parameters for the lower bound: $e=f=g=h=i=j=k=m=n=o=q=0$ and $a_1=a_2=a_3=a_4=a=b=c=d=l=p=r=s=t=3/7$. This implies\smallskip

$Q\left(3/7,3/7,3/7,3/7\right)=Q\left(3/7,3/7,3/7,6/7\right)=Q\left(3/7,3/7,6/7,3/7\right)=Q\left(3/7,3/7,6/7,6/7\right)=$

$Q\left(3/7,6/7,3/7,3/7\right)=Q\left(3/7,6/7,3/7,6/7\right)=Q\left(3/7,6/7,6/7,3/7\right)=Q\left(6/7,3/7,3/7,3/7\right)=$

$Q\left(6/7,3/7,3/7,6/7\right)=Q\left(6/7,3/7,6/7,3/7\right)=Q\left(6/7,6/7,3/7,3/7\right)=0$\smallskip

\noindent
and\smallskip

$Q\left(3/7,6/7,6/7,6/7\right)=Q\left(6/7,3/7,6/7,6/7\right)=Q\left(6/7,6/7,3/7,6/7\right)=Q\left(6/7,6/7,6/7,3/7\right)=$

$Q\left(6/7,6/7,6/7,6/7\right)=3/7$,\smallskip

\noindent
whence the 4-box $B_1=[3/7,6/7]^4$ is obtained, which is symmetrical with respect to the main diagonal, as in the bi- and trivariate cases.

On the other hand, the solution for the parameters obtained for the upper bound is $e=f=g=i=m=0$, $a_1=a_2=a_3=a_4=a=b=c=d=h=j=k=l=n=o=p=q=r=s=1/2$ and $t=1$, whence the 4-box is given by $B_2=\left[1/2,1\right]^4$ and the values on the vertices of $B_2$ are given by:\smallskip

$Q(1/2,1/2,1,1)=Q(1,1,1/2,1/2)=Q(1/2,1,1/2,1)=Q(1/2,1,1,1/2)=Q(1,1/2,1/2,1)=$

$Q(1,1/2,1,1/2)=1/2$.\smallskip

\noindent
Obviously we have\smallskip

$Q(1/2,1,1,1)=Q(1,1/2,1,1)=Q(1,1,1/2,1)=Q(1,1,1,1/2)=1/2$ and $Q(1,1,1,1)=1$.
\smallskip

Also, we want to stress that the software does not indicate a unique solution for the lower bound, nor that there are infinitely many solutions for the upper bound, as occurs in the trivariate case, although it is probable that a similar result can be obtained for the tetravariate case, remaining this as an open problem.

We now provide two examples of 4-quasi-copulas, $Q_1$ and $Q_2$, for which $V_{Q_{1}}\left([3/7,6/7]^4\right)=-9/7$ and $V_{Q_{2}}\left([1/2,1]^4\right)=2$.

\begin{example}\label{ex:ex1}Let $Q_1$ be the 4-quasi-copula whose mass is spread uniformly on $[0,1]^4$ in the following manner: $3/7$ of (positive) mass on the 4-boxes $[3/7,6/7]^3\times[0,3/7]$, $[3/7,6/7]^2\times[0,3/7]\times[3/7,6/7]$, $[3/7,6/7]\times[0,3/7]\times[3/7,6/7]^2$ and $[0,3/7]\times[3/7,6/7]^3$; $1/7$ of (positive) mass on the 4-boxes $[3/7,6/7]^3\times[6/7,1]$, $[3/7,6/7]^2\times[6/7,1]\times[3/7,6/7]$, $[3/7,6/7]\times[6/7,1]\times[3/7,6/7]^2$ and $[6/7,1]\times[3/7,6/7]^3$; $-9/7$ of (negative) mass on the 4-box $[3/7,6/7]^4$; and $0$ on the remaining 4-boxes. In order to visualize, in some way, the 4-quasi-copula $Q_1$, we show the mass distribution of one of its identical trivariate margins, which are 3-quasi-copulas, in Figure \ref{fig:Q1}.
\begin{figure}
\hspace{1.5cm}\epsfig{file=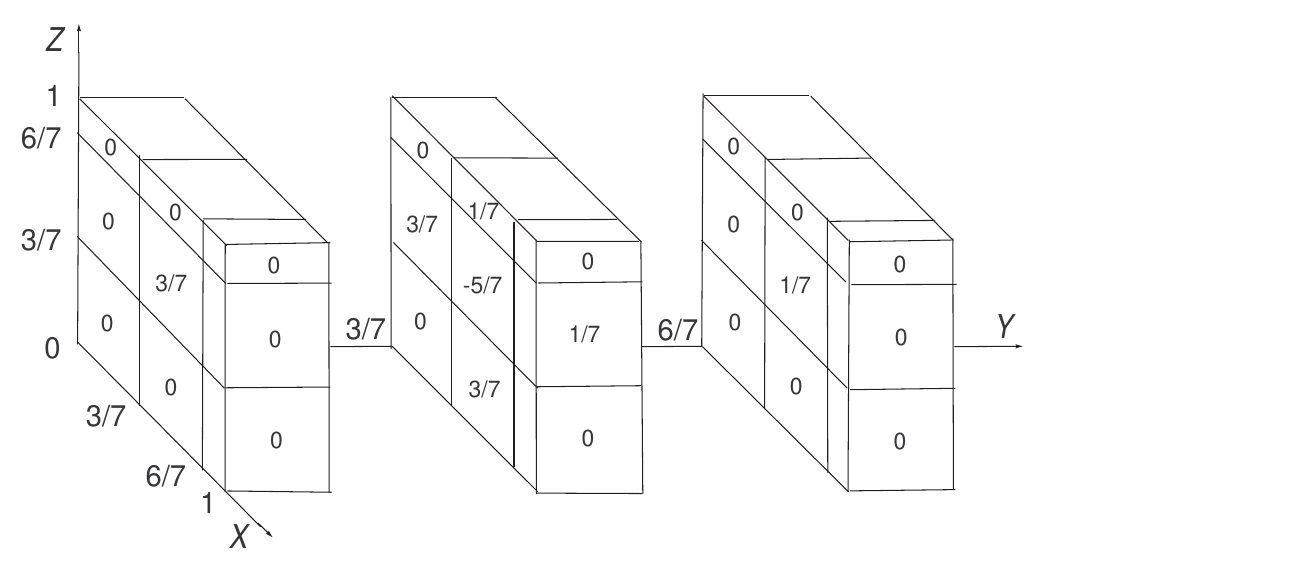,width=\textwidth}
\caption{Mass distribution of the trivariate marginal of the 4-quasi-copula $Q_1$ in Example \ref{ex:ex1}.}\label{fig:Q1}
\end{figure}
\end{example}

\begin{example}Let $Q_2$ be the 4-quasi-copula whose mass is spread uniformly on $[0,1]^4$ in the following manner: $2$ of (positive) mass on the 4-box $[1/2,1]^4$; $1/2$ of (positive) mass on the 4-boxes $[0,1/2]^2\times[1/2,1]^2$, $[0,1/2]\times[1/2,1]\times[0,1/2]\times[1/2,1]$, $[0,1/2]\times[1/2,1]^2\times[0,1/2]$, $[1/2,1]\times[0,1/2]^2\times[1/2,1]$, $[1/2,1]\times[0,1/2]\times[1/2,1]\times[0,1/2]$ and $[1/2,1]^2\times[0,1/2]^2$; $-1$ of (negative) mass on the 4-boxes $[0,1/2]\times[1/2,1]^3$, $[1/2,1]\times[0,1/2]\times[1/2,1]^2$, $[1/2,1]^2\times[0,1/2]\times[1/2,1]$ and $[1/2,1]^3\times[0,1/2]$; and $0$ on the remaining 4-boxes. We note that the (identical) trivariate margins of the 4-quasi-copula $Q_2$ correspond to a 3-quasi-copula $Q_3$ for which $V_{Q_3}\left([1/2,1]^3\right)=1$: see \cite[Example 2]{DeBaets2007}.
\end{example}

\section{Conclusion}\label{sec:discus}

$\phantom{999}$In this note, we have shown that, ``surprisingly'', the extreme values of the mass distribution in a 4-box of a 4-quasi-copula are $-9/7$ (less than $-1$) and 2 (more than 1), showing the existing differences with respect to the bi- and trivariate cases. In view of the results obtained in all these cases, a natural question arises: what are the extremes of the mass distribution associated with an $n$-quasi-copula for $n>4$? \cite[Open Problem 5]{Arias2020}. The answer seems not to be simple, even more considering that the use of linear programming does not seem feasible. In any case, we conjecture, at least, for the lower bound, the following: Let $n$ be a natural number such that $n\ge 2$, let $Q$ be an $n$-quasi-copula, and let $B$ be an $n$-box in $[0,1]^n$. Then
$$V_Q(B)\ge -\frac{(n-1)^2}{2n-1};$$
moreover, the equality implies
$$B=\left[\frac{n-1}{2n-1},\frac{2n-2}{2n-1}\right]^n.$$
\bigskip\smallskip

\noindent{\large\bf Acknowledgements}
\medskip

The author acknowledges the support of the program FEDER-Andaluc\'ia 2014-2020 under research project UAL2020-AGR-B1783 and project PID2021-122657OB-I00 by the Ministerio de Ciencia e Innovaci\'on (Spain).

\end{document}